\documentclass[12pt]{article}
\usepackage{amsmath, amssymb}

\begin{document}

\title{\bf A remark about positive polynomials. }
\author{Olga M. Katkova, Anna M. Vishnyakova\\
Dept. of Math., Kharkov National University, \\
Svobody sq., 4, 61077, Kharkov, Ukraine \\
e-mail: olga.m.katkova@univer.kharkov.ua, \\
anna.m.vishnyakova@univer.kharkov.ua } \maketitle

\begin{abstract}
The following theorem is proved.

{\bf Theorem.} {\it Let $P(x) = \sum_{k=0}^{2n} a_k x^k$ be a
polynomial with positive coefficients. If the inequalities
$\frac{a_{2k+1}^2}{a_{2k}a_{2k+ 2}} <
\frac{1}{cos^2(\frac{\pi}{n+2})} $ hold for  all $ k=0, 1, \ldots,
n-1, $ then $P(x)>0$ for every $x\in\mathbb{R} $  .}

We  show that the constant $\frac{1}{cos^2(\frac{\pi}{n+2})}$ in
this theorem could not be increased. We also present some
corollaries of this theorem.

\end{abstract}

{\it 2000 Mathematics Subject Classification  30C15; 26C10.}

\section{Introduction and statement of results.}

Positive polynomials arise in many important branches of
mathematics. In this note we give a simple sufficient condition
for an even degree polynomial with positive coefficients to be
positive on the real line. Before we formulate the main theorem we
will mention two results which have in some sense similar
character.

In 1926, Hutchinson \cite[p.327]{hut} extended the work of
Petrovitch \cite{pet} and Hardy \cite{har1} or \cite[pp.
95-100]{har2} and proved the following theorem.

{\bf Theorem A.} { \it Let $P(x) = \sum_{k=0}^{n} a_k x^k$ be a
polynomial with positive coefficients. If the inequalities
\begin{equation}
\frac{a_k^2}{a_{k-1}a_{k+1}}\geq 4 ,\ k=1, 2, \ldots, n-1,
\end{equation}
hold then all zeros of $P(x)$ are real. }

In \cite{kurtz} it was proved that the constant $4$ in Theorem A
is sharp.

In \cite{kavi} the authors of this note have found the smallest
possible constant $d_n> 0$ such that if coefficients of
$P(x)=\sum_{k=0}^{n} a_k x^k$ are positive and satisfy the
inequalities $ \frac{a_k^2}{a_{k-1}a_{k+1}} > d_n,
 \ k=1, 2, \ldots, n-1 , $ then
$P(x)$ is Hurwitz. We remind that a real polynomial is called
Hurwitz (stable) if all its zeros have negative real parts.

The following theorem is the main result of this work.

{\bf Theorem 1.} {\it Let $P(x) = \sum_{k=0}^{2n} a_k x^k$ be a
polynomial with positive coefficients. If the inequalities
$$\frac{a_{2k+1}^2}{a_{2k}a_{2k+ 2}} <
\frac{1}{cos^2(\frac{\pi}{n+2})} $$ hold for  all $ k=0, 1,
\ldots, n-1, $ then $P(x)>0$ for every $x\in\mathbb{R} $  .}

The following  Theorem shows that the constant
$\frac{1}{cos^2(\frac{\pi}{n+2})}$ in Theorem 1 is sharp for every
$n\in\mathbb{N}$.

{\bf Theorem 2.} {\it For every $n\in\mathbb{N}$ there exists a
polynomial $Q(x)= \sum_{k=0}^{2n} a_k x^k$  with positive
coefficients under conditions $$\frac{a_{2k+1}^2}{a_{2k}a_{2k+ 2}}
= \frac{1}{cos^2(\frac{\pi}{n+2})}, \quad   k=0, 1, \ldots, n-1,
$$ and the polynomial  $Q(x)$ has not less than two real
zeros. }

The following statement is a simple corollary of Theorem 1.

{\bf Corollary 1.} {\it Let $P(x) = \sum_{k=0}^{2n+1} a_k x^k$ be
a polynomial with positive coefficients. If the inequalities
$$\frac{a_{2k}^2}{a_{2k-1}a_{2k+ 1}} < \frac{4k^2-1}{4k^2}\cdot
\frac{1}{cos^2(\frac{\pi}{n+2})} $$ hold for  all $ k= 1, 2,
\ldots, n, $ then $P(x)$ has only one real zero (counting
multiplicities). }

We show that the constants in the last statement is also sharp for
every $n\in\mathbb{N}$.

{\bf Theorem 3.} {\it For every $n\in\mathbb{N}$ there exists a
polynomial $Q(x)= \sum_{k=0}^{2n+1} a_k x^k$  with positive
coefficients under conditions $$\frac{a_{2k}^2}{a_{2k-1}a_{2k+ 1}}
= \frac{4k^2-1}{4k^2}\cdot \frac{1}{cos^2(\frac{\pi}{n+2})}, \quad
k= 1, 2, \ldots, n, $$ and the polynomial $Q(x)$ has not less than
three real zeros. }

\section {Proof of Theorem 1.}

Let $P(x) = \sum_{k=0}^{2n} a_k x^k$ be a polynomial with positive
coefficients. Let us consider a quadratic form
\begin{eqnarray}
\label{f2} & Q_P(x_0, x_1, x_2,  \ldots , x_n)= a_0 x_0^2 +a_1
x_0x_1 + a_2x_1^2+a_3x_1x_2+a_4x_2^2+
 \ldots \\ \nonumber &
 +a_{2n-2}x_{n-1}^2 +a_{2n-1}x_{n-1}x_n +a_{2n}x_n^2 = \sum_{k=0}^n a_{2k}x_k^2
 + \sum_{k=0}^{n-1}a_{2k+1}x_kx_{k+1} .
\end{eqnarray}
For every $x\in\mathbb{R}$ we have
\begin{equation}
\label{f3} P(x) = Q_P(1, x, x^2 \ldots , x^n).
\end{equation}
Thus, if the quadratic form $Q_P$ is positive-definite, then for
every $x\in\mathbb{R}$ we have $P(x)>0$. It remains to prove that
under assumptions of Theorem 1 the quadratic form $Q_P$ is
positive-definite. A following $(n+1) \times (n+1)$ matrix
corresponds to the quadratic form $Q_P$

\begin {equation}
M_{Q_p}:= \left\|
  \begin{array}{ccccccc}
   a_0 & \frac{a_1}{2} & 0 & 0 &\ldots &0&0\\
   \frac{a_1}{2}   & a_2 & \frac{a_3}{2} & 0 &\ldots&0&0 \\
   0   &  \frac{a_3}{2}  & a_4 & \frac{a_5}{2} &0&\ldots&0 \\
   \vdots&\vdots&\vdots&\vdots&\ldots&\vdots&\vdots\\
   0&0&\ldots&0& \frac{a_{2n-3}}{2}   & a_{2n-2} & \frac{a_{2n-1}}{2}\\
   0  &  0  &0 &\ldots &0& \frac{a_{2n-1}}{2} & a_{2n} \\
     \end{array}
 \right\|.
\label{f4}
\end {equation}

By Sylvester's Criterion for positive definiteness we need to show
that all leading principal minors of the matrix $M_{Q_p}$ are
positive. To do this we will use the following theorem from
\cite{kv}.

{\bf Theorem A.} {\it Let  $M = (a_{ij})$ be an  $m \times m$
matrix with the properties \\
(a) $ a_{ij} > 0 \  (1 \leq i,j \leq m)$ and \\
(b) $a_{ij}a_{i+1,j+1} > 4 \cos ^2 \frac {\pi}{m+1}\
a_{i,j+1}a_{i+1,j} \
(1 \leq i,j \leq m-1) .$\\
Then all minors of $M$ are positive. }

In \cite{kv} it is also shown that the constant $c_m := 4 \cos ^2
\frac {\pi}{m+1}$ in the statement of Theorem A is the smallest
possible not only in the class of $m \times m$ matrices with
positive entries but in the classes of $m \times m$ Toeplitz
matrices and of $m \times m$ Hankel matrices.

Consider the following $(n+1)\times (n+1)$ symmetrical  Toeplitz
matrix
\begin {eqnarray}
& T( \varepsilon_1,\ldots , \varepsilon_{n-1}):=  \left\|
  \begin{array}{ccccccc}
   a_0 & \frac{a_1}{2} & \varepsilon_1 & \varepsilon_2 &\ldots &\varepsilon_{n-2}
   &\varepsilon_{n-1}\\
   \frac{a_1}{2}   & a_2 & \frac{a_3}{2} & \varepsilon_1 &\ldots&
   \varepsilon_{n-3}&\varepsilon_{n-2} \\
   \varepsilon_1   &  \frac{a_3}{2}  & a_4&\frac{a_5}{2}&
   \ldots&\varepsilon_{n-4}&
   \varepsilon_{n-3} \\
   \vdots&\vdots&\vdots&\vdots&\ldots&\vdots&\vdots\\
   \varepsilon_{n-2}&\varepsilon_{n-3}&\ldots&
   \varepsilon_1&\frac{a_{2n-3}}{2}&
    a_{2n-2} & \frac{a_{2n-1}}{2}\\
   \varepsilon_{n-1}  &  \varepsilon_{n-2}  &\varepsilon_{n-3} &
   \ldots&\varepsilon_1& \frac{a_{2n-1}}{2} & a_{2n} \\
     \end{array}
 \right\|,
\label{f5}
\end {eqnarray}
where  $\varepsilon_1 > \varepsilon_2 > \cdots >
\varepsilon_{n-1}>0$ will be chosen in such a way that the matrix
$T( \varepsilon_1,\ldots , \varepsilon_{n-1})$ will satisfy the
assumptions of Theorem A. At first we will choose $\varepsilon_1$
such that
$$ \frac{a_{2j-1}a_{2j+1}}{4} >  4 \cos ^2
\frac {\pi}{n+2}\ a_{2j} \varepsilon_1 ,\quad  j=1, 2, \ldots ,
n-1 .$$
After that we will choose $\varepsilon_2$ such that
$$ \varepsilon_1^2 > 4 \cos ^2 \frac {\pi}{n+2}\ a_{2j+1}
\varepsilon_2 ,\quad  j=1, 2, \ldots , n-2 .$$
Then we will choose
$\varepsilon_3 > \varepsilon_4 > \cdots > \varepsilon_{n-1}>0$ one
after another such that
$$ \varepsilon_{j}^2 > 4 \cos ^2\frac {\pi}{n+2}\varepsilon_{j-1}\varepsilon_{j+1} ,
\quad  j=2, 3, \ldots , n-2.$$ For our choice of $\varepsilon_1>
\varepsilon_2 > \cdots > \varepsilon_{n-1}>0$ and under the
assumptions on $a_0, a_1, \ldots , a_{2n}$ the matrix $T(
\varepsilon_1,\ldots , \varepsilon_{n-1})$ satisfies the
assumptions of Theorem A. So by Theorem A all minors of $T(
\varepsilon_1,\ldots , \varepsilon_{n-1})$ are positive. Tending
$\varepsilon_{n-1} \to 0,$  $\varepsilon_{n-2}\to 0, \ldots ,$
$\varepsilon_{1}\to 0$ we obtain that all minors of $M_{Q_p}$ are
nonnegative. It remains to prove that all leading principal minors
of the matrix $M_{Q_p}$
$$\Delta_1(M_{Q_p}) = a_0,\    \Delta_2 (M_{Q_p}) =
\det \left( \begin{array}{cc}
   a_0 & \frac{a_1}{2} \\
   \frac{a_1}{2}   & a_2  \\
     \end{array}  \right) ,
\ldots , \   \Delta_{n+1}(M_{Q_p})= \det M_{Q_p}  $$ are positive.

Suppose there is a leading principal minor of  $M_{Q_p}$ which is
equal to zero. Denote by $j$ the minimal number of a leading
principal minor which is equal to zero. Since $\Delta_1(M_{Q_p}) =
a_0 >0$, we have $j\geq 2$ and $\Delta_{j-1}(M_{Q_p})>0,
\Delta_j(M_{Q_p})=0.$ Let us consider a polynomial
$P_{\varepsilon} (x) = P(x) -\varepsilon x^{2j-2}$ where
$\varepsilon > 0$ is so small that $P_{\varepsilon} (x)$ satisfies
the assumptions of Theorem 1. As we have proved it implies that
all minors of a corresponding matrix $ M_{Q_{p_{\varepsilon}}}$
are nonnegative, in particular
$$ \Delta_j(M_{Q_{p_{\varepsilon}}})= \det \left(
  \begin{array}{ccccccc}
   a_0 & \frac{a_1}{2} & 0 & 0 &\ldots &0&0\\
   \frac{a_1}{2}   & a_2 & \frac{a_3}{2} & 0 &\ldots&0&0 \\
   0   &  \frac{a_3}{2}  & a_4 & \frac{a_5}{2} &0&\ldots&0 \\
   \vdots&\vdots&\vdots&\vdots&\ldots&\vdots&\vdots\\
   0&0&\ldots&0& \frac{a_{2j-5}}{2}   & a_{2j-4} & \frac{a_{2j-3}}{2}\\
   0  &  0  &0 &\ldots &0& \frac{a_{2j-3}}{2} & a_{2j-2} -\varepsilon \\
     \end{array}
 \right) \geq 0 .$$
Since $\Delta_{j-1}(M_{Q_p})>0$ we conclude that the determinant $
\Delta_j(M_{Q_{p_{\varepsilon}}})$ is strictly decreasing in
$\varepsilon > 0$ and so $\Delta_j(M_{Q_p})$ (which is equal to
$\Delta_j(M_{Q_{p_{\varepsilon}}})$ for $\varepsilon = 0$) is
strictly positive. Thus there are no zero leading principal minors
of $M_{Q_p}$, all leading principal minors of  $M_{Q_p}$ are
positive. By Sylvester's Criterion it means that the quadratic
form $Q_P$ is positive-definite, and in particular from (\ref{f3})
we obtain that $P(x)>0$ for every $x\in\mathbb{R}$.

Theorem 1 is proved.  $\Box$

\section {Proof of Theorems 2 and 3.}

{\it Proof of Theorem 2.}  Let us fix an arbitrary $n\in
\mathbb{N}$  and denote by $\alpha :=\frac{\pi}{n+2}.$ Consider
the following polynomial
$$ Q(x):= \sum_{k=1}^n \sin k\alpha \sin(k+1)\alpha \
(1+x)^2x^{2k-2}.$$ Obviously $Q(x)$ is a polynomial with positive
coefficients of degree $2n$ and $-1$ is a root of $Q$ of
multiplicity not less than $2$. We have
\begin{eqnarray}
\label{f6} & Q(x)= \sum_{k=1}^n \sin k\alpha \sin(k+1)\alpha \
(x^{2k-2} +x^{2k})  \\  \nonumber & + 2 \sum_{k=1}^n \sin k\alpha
\sin(k+1)\alpha \   x^{2k-1}  =
\sum_{k=0}^{n-1} \sin (k+1)\alpha \sin(k+2)\alpha \   x^{2k} \\
\nonumber & + \sum_{k=1}^n \sin k\alpha \sin(k+1)\alpha \  x^{2k}
 + 2 \sum_{k=1}^n \sin k\alpha \sin(k+1)\alpha \
x^{2k-1} \\  \nonumber & = \sin\alpha \sin 2\alpha +
\sum_{k=1}^{n-1} \sin(k+1)\alpha ( \sin(k+2)\alpha + \sin k\alpha) x^{2k} + \\
\nonumber & \sin (n+1)\alpha (\sin n\alpha +\sin(n+2)\alpha)
x^{2n} + 2 \sum_{k=1}^n \sin k\alpha \sin(k+1)\alpha \   x^{2k-1}
\\  \nonumber & = 2\sin^2\alpha \cos \alpha +
2\sum_{k=1}^{n-1}\sin^2(k+1)\alpha \cos\alpha \   x^{2k} +
2\sin^2(n+1)\alpha \cos\alpha \   x^{2n}+\\  \nonumber &  2
\sum_{k=1}^n \sin k\alpha \sin(k+1)\alpha \   x^{2k-1} = \\
\nonumber & 2\sum_{k=0}^{n}\sin^2(k+1)\alpha \cos\alpha \
x^{2k}+2 \sum_{k=1}^n \sin k\alpha \sin(k+1)\alpha \   x^{2k-1}
\end{eqnarray}
(we use the fact that $\sin(n+2)\alpha=0).$ So if we define by
$a_j, \ j = 0, 1,  \ldots , 2n ,$ the coefficients of $Q$ then
$$a_{2k}= 2 \sin^2(k+1)\alpha \cos\alpha ,
\ a_{2k-1}= 2 \sin k\alpha \sin(k+1)\alpha .$$ We have

$$\frac{a_{2k+1}^2}{a_{2k}a_{2k+ 2}}= \frac{4 \sin^2 (k+1)\alpha
\sin^2(k+2)\alpha}{2 \sin^2(k+1)\alpha \cos\alpha \cdot 2
\sin^2(k+2)\alpha \cos\alpha} = \frac{1}{\cos^2\alpha}$$ for
$k=0,1, \ldots, n-1.$

Theorem 2 is proved.  $\Box$

{\it Proof of Theorem 3.} Let us fix an arbitrary $n\in
\mathbb{N}$  and denote by $\alpha :=\frac{\pi}{n+2}.$ Let us
consider the following primitive function for the polynomial
$Q(x)$ constructed in the proof of Theorem 2:
$$ H(x)= \sum_{k=1}^n \sin k\alpha \sin(k+1)\alpha
\left(\frac{x^{2k-1}}{2k-1}+ 2\frac{x^{2k}}{2k}
+\frac{x^{2k+1}}{2k+1}\right) .$$ We have
$$ H(-1) = \sum_{k=1}^n \sin k\alpha \sin(k+1)\alpha
\left(-\frac{1}{2k-1}+ \frac{1}{k} - \frac{1}{2k+1}\right)$$
$$ = \sum_{k=1}^n \sin k\alpha \sin(k+1)\alpha \left(
\frac{-1}{k(2k-1)(2k+1)}\right) <0 .$$

So the following polynomial
$$ S(x) = H(x) - H(-1)$$
is a polynomial with positive coefficients of degree $2n+1$ and
$-1$ is a root of $S$ of multiplicity not less than $3$. Using
(\ref{f6}) we can rewrite $S(x)$ in the form
$$ S(x)= - H(-1) +2\sum_{k=0}^{n}\sin^2(k+1)\alpha \cos\alpha \
\frac{x^{2k+1}}{2k+1}+2 \sum_{k=1}^n \sin k\alpha \sin(k+1)\alpha
\ \frac{x^{2k}}{2k} .$$ So if we define by $b_j, \ j = 0, 1,
\ldots , 2n+1,$ the coefficients of $Q$ then $b_0=- H(-1)$ and
$$b_{2k+1}= \frac{2 \sin^2(k+1)\alpha \cos\alpha}{2k+1} ,\   k= 0, 1, \ldots ,
n ;$$ $$ b_{2k}= \frac{2 \sin k\alpha \sin(k+1)\alpha}{2k}, \   k=
1, 2 \ldots , n.$$ We have
$$\frac{b_{2k}^2}{b_{2k-1}b_{2k+ 1}}
= \frac{4 \sin^2 k\alpha \sin^2(k+1)\alpha}{4k^2}\cdot
\frac{2k-1}{2 \sin^2 k\alpha \cos\alpha} \frac{2k+1}{2
\sin^2(k+1)\alpha \cos\alpha}$$
$$ = \frac{4k^2-1}{4k^2}\cdot \frac{1}{\cos^2\alpha}, \    k=
1, 2 \ldots , n.$$

Theorem 3 is proved.  $\Box$

\end{document}